\documentclass[12pt]{amsart}
\begingroup
\swapnumbers
\theoremstyle{plain}

\newtheorem{thm}{Theorem}[section]
\newtheorem{cor}[thm]{Corollary}
\newtheorem{prop}[thm]{Proposition}
\newtheorem{lem}[thm]{Lemma}

\theoremstyle{definition}
\newtheorem{dfn}[thm]{Definition}
\newtheorem{rmk}[thm]{Remark}
\newtheorem{rmks}[thm]{Remarks}
\newtheorem{egz}[thm]{Example}

\newtheorem{gen}{}
\endgroup

\numberwithin{equation}{section}

\newcommand {\Bh}{\mathcal{B}(\mathcal{H})}  
\newcommand {\Bk}{\mathcal{B}(\mathcal{K})}  
\newcommand {\Bhk}{\mathcal{B}(\mathcal{H},\mathcal{K})}  
\newcommand {\Bkoh}{\mathcal{B}(\mathcal{K}\oplus\mathcal{H})}  
\newcommand {\B}{\mathcal{B}}  
\newcommand {\HH}{\mathcal{H}}  
\newcommand {\KK}{\mathcal{K}} 
\newcommand {\ld}{\lambda} 

\newcommand {\N}{\mathbb{Z}_{+}} 
\newcommand {\R}{\mathbb{R}} 
\newcommand {\C}{\mathbb{C}} 

\newcommand {\om}{\omega}
\newcommand {\som}{S_{\omega}}

\def\firsta{\  {\text {\rm (a)}}\ \ }
\def\conda{\smallskip  \noindent {\text {\rm (a)}}\ \ }
\def\condb{\smallskip  \noindent {\text {\rm (b)}}\ \ }
\def\condc{\smallskip  \noindent {\text {\rm (c)}}\ \ }

\newcommand {\be}{\begin{equation}}
\newcommand {\ee}{\end{equation}}
\newcommand {\beq}{\begin{eqnarray*}}
\newcommand {\eeq}{\end{eqnarray*}}

\begin{document}

\title[Perturbations of operators and the commutator equation]
{Perturbations of operators similar to contractions and the commutator equation}

\author{C. Badea}
\address{D\'epartement de Math\'ematiques, UMR 8524 au CNRS, 
Universit\'e de Lille 1,
F--59655 Villeneuve d'Ascq, France}
\email{Catalin.Badea@math.univ-lille1.fr}
\urladdr{http://math.univ-lille1.fr/\~{}badea/}
\date{}

\begin{abstract}  
Let $T$ and $V$ be two Hilbert space contractions and let $X$ be a linear bounded operator. 
It was proved by C.~Foias and J.P.~Williams that in certain
cases the operator block matrix $R(X;T,V)$ (equation
(\ref{eq:rx}) below) is similar to a contraction if and only if the
commutator equation $X = TZ-ZV$ has a bounded solution $Z$. We characterize here the similarity 
to contractions of some operator matrices $R(X;T,V)$ in terms of growth conditions or of 
perturbations of $R(0;T,V)=T\oplus V$.
\end{abstract}

\maketitle


\section{Introduction}
\thispagestyle{empty}
A bounded linear operator $T \in \Bh$, acting on a Hilbert space $\HH$, is said to 
be \emph{polynomially bounded} if there exists a constant $M$ such that the
inequality 
$$
\|p(T)\| \le M\|p\|_{\infty} =: M \sup\{ |p(z)| : |z| = 1\}
$$
holds for all polynomials $p \in \C[z]$. It is said to be  
\emph{power bounded} if the same inequality holds for all monomials 
$p_n(z) = z^n$. 

The von Neumann
inequality implies that every operator $T$ similar to a contraction is polynomially
bounded. Recall that $T$ is said to be \emph{similar to a contraction} if there exists an invertible
operator $L\in \Bh$ such that $\|L^{-1}TL\| \leq 1$.

There are power bounded operators which are not polynomially bounded 
(Foguel \cite{foguel} ; see
also \cite{davie, peller2,bozejko}) as well as
polynomially bounded operators which are not similar to a contraction 
(Pisier \cite{pisier:jams} ; see also \cite{DaPa}). 
Both Foguel's and Pisier's counterexamples are operators of the following type~:
$$
R(X;V^*,V) = \left[\begin{array}{cc}
            V^* & X  \\
            0 & V
        \end{array}\right] \in
        \mathcal{B}(\mathcal{K}\oplus\mathcal{K}),$$
where $V\in \Bk$ is a pure isometry (i.e. a unilateral shift), $V^*$ is its adjoint and $X$ is an
operator in $\Bk$. In Foguel's counterexample, $\KK = \ell^2$ and $X$ is a 
suitable diagonal projection onto a subspace of $\ell^2$, while in Pisier's counterexample, $\KK =
\ell^2(\HH)$, $V$ is the unilateral shift of infinite multiplicity $\dim \HH$ and $X$ is a suitable
Hankel operator with suitable operator-valued entries. 

Let $T \in \Bk$, $V\in \Bh$ and $X\in \Bhk$ be three Hilbert space operators. In this paper we will 
consider operators of the form 
\begin{equation}\label{eq:rx}
R(X) = R(X;T,V) = \left[\begin{array}{cc}
            T & X  \\
            0 & V
        \end{array}\right] \in \mathcal{B}(\mathcal{K}\oplus\mathcal{H}).
\end{equation}
Suppose that $T$ is a coisometry (i.e. its adjoint $T^*$ is an isometry) 
and $V$ is a contraction or that $T$ is a
contraction and $V$ is an isometry.
It was proved by C. Foias and J.P. Williams \cite{foias/williams} (cf. \cite{C+3,paulsen:yoneda,clark}) 
that in these cases
$R(X)$ is similar to a contraction if and only if the 
(generalized) commutator equation 
\begin{equation}\label{eq:comm}
X = TZ -ZV
\end{equation}
has a bounded solution $Z \in \Bhk$. This implies that $R(X)$ is similar to a contraction if
and only if $R(X)$ is similar to $R(0)=T\oplus V$, if and only if the commutator equation 
has a bounded solution. This follows from the matrix
identity
$$\left[\begin{array}{cc}
            I & -Z  \\
            0 & I
        \end{array}\right] \left[\begin{array}{cc}
            T & 0  \\
            0 & V
        \end{array}\right]\left[\begin{array}{cc}
            I & Z \\
            0 & I
        \end{array}\right] =
        \left[\begin{array}{cc}
            T & TZ-ZV  \\
            0 & V
        \end{array}\right] .$$
Here $I$ denotes the identity operator, possibly on different spaces. 

Note that, given $T$ and $V$, 
the equivalence between the similarity of $R(X)$ to $R(0)$ and the 
solvability of
the commutator equation (\ref{eq:comm}) holds for all operators 
$T$ and $V$ if $\HH$ and $\KK$ are finite dimensional
(Roth's theorem \cite{roth}). It also holds for some, but not all,  
pairs $(T,V)$ in the infinite dimensional situation 
\cite{rosenblum},\cite{bhatia/rosenthal}.
We also note that the equation (\ref{eq:comm}) is 
sometimes called in the literature
the Sylvester equation, while the operator $Z\to TZ-ZV$ is concurrently called the Rosenblum
operator, an elementary operator, 
a generalized commutator or a generalized derivation. 
We refer to the survey paper \cite{bhatia/rosenthal} 
by Bhatia and Rosenthal for more information 
concerning the commutator equation (\ref{eq:comm}).  

The aim of the present paper is to give different characterizations of operators $R(X)$ 
similar to contractions. In all cases the Foias-Williams' theorem will be applicable, 
and so these results can be viewed as criteria for the solvability of the
commutator equation. All characterizations will be applicable to the case when $T$ is a 
coisometry and $V$ is a isometry, not necessarily acting on the same space.

In the first such result (Theorem \ref{thm:gr}), $T$ is supposed to be a right
invertible contraction and $V$ an isometry or $T$ a coisometry and $V$ a left invertible 
contraction. The result says that in both cases $R(X)$ is similar to a contraction
if and only if $X$ can be decomposed as the sum of two bounded operators satisfying growth
conditions or product identities. 

A different characterization for operators $R(X;T,V)$, with 
$T$ a coisometry and $V$ a weighted
unilateral shift of arbitrary multiplicity, 
is given in Theorem \ref{thm:main3}.
This time the necessary and sufficient condition is that 
$R(X)$ is power bounded and it can be written 
as a zero-product perturbation
by a nilpotent 
of an operator which is near (in a certain sense) to $R(0)$. 
Basically the same characterization holds (Theorem \ref{thm:main1}) if $T$ is a right
invertible contraction and $V$ is an isometry.
These results are related to a previous
result of Apostol \cite{apostol:pams} concerning the commutator equation. 
Note also that other variations are possible.

Starting with work of K.O. Friedrichs (cf.
\cite[Part III, XX.2.2]{dunford/schwartz}), several results of the type  "a perturbation of an operator $C$ 
is similar to $C$" exist in the literature.
The characterizations described above can be interpreted as results of the following type~:
"perturbations of operators near $R(0)= T\oplus V$ are similar to $R(0)$".

The paper is organized as follows. In the next section we study the commutator equation and 
Theorem \ref{thm:gr} is proved. The zero-product perturbations are introduced in section 3 and 
the stability of the class of operators similar to contractions under such perturbations is 
studied.
In section 4 we introduce the notion of $\beta$-quadratically near operators modulo subspaces, 
extending \cite[Definition 2.5]{b:qn}.
This notion is used in the main results of section 5.
\section{Growth conditions and the commutator equation}
The proofs of the following result, and that of Theorem \ref{thm:blbis} below, 
are similar to the proof of \cite[Theorem 4.1]{bp}. We refer to 
this paper for more references and to \cite[p. 9]{bhatia/rosenthal} 
for the significance of the condition (\ref{eq:gr}). 
\begin{thm}\label{thm:bl}
Let $V\in \Bh$ be an isometry, let $T\in \Bk$ be a power bounded operator and let 
$X\in \Bhk$.
Then the commutator equation
$$X = TZ -ZV$$
has a bounded solution $Z\in \Bhk$ satisfying $Z(I - VV^*) = 0$ if and only if 
\begin{equation}\label{eq:gr}
\sup_{n}\|\sum_{j=0}^n T^jXV^{*j+1}\| < \infty .
\end{equation}
\end{thm}
\begin{proof}
Suppose the equation $X = TZ -ZV$ has a bounded solution $Z\in \Bhk$ satisfying
$Z = ZVV^*$.
Then 
\beq
\sum_{j=0}^n T^jXV^{*j+1} & = & \sum_{j=0}^n T^j(TZ - ZV)V^{*j+1} \\
& = & \sum_{j=0}^n T^{j+1}ZV^{*j+1}  - \sum_{j=0}^n T^{j}ZVV^*V^{*j} \\
& = & \sum_{j=0}^n T^{j+1}ZV^{*j+1}  - \sum_{j=0}^n T^{j}ZV^{*j} \\
 & = & T^{n+1}ZV^{*n+1} - Z .
\eeq
Therefore, $\|T^{n+1}\| \leq M$ implies
$$
\sup_{n}\|\sum_{j=0}^n T^jXV^{*j+1}\| \leq (M+1)\|Z\| .
$$

Suppose now that (\ref{eq:gr}) holds. Let $\mathcal{L}$ be a
Banach limit \cite[Part I,p. 73]{dunford/schwartz}, that is a bounded
linear functional on $\ell_{\infty}(\C)$ such that ${\bf 1} =
\mathcal{L}({\bf 1}) = \|\mathcal{L}\|$ and $\mathcal{L}((x_{n+1})_{n\geq
0}) = \mathcal{L}((x_{n})_{n\geq 0})$ for every $(x_{n})_{n\geq 0} \in
\ell_{\infty}(\C)$. Here ${\bf 1} = (1,1, \ldots)$.

Consider the linear operator $Z : \mathcal{H} \to \mathcal{K}$ given by
$$\langle Zh , k\rangle = - \mathcal{L}\left( \langle
\sum_{j=0}^n T^jXV^{*j+1} h , k\rangle\right) .$$
Then (\ref{eq:gr}) shows that $Z$ is well-defined and bounded.

We have 
\beq
\langle (TZ - ZV)h,k\rangle & = & \langle
    Zh,T^*k\rangle - \langle ZVh,k\rangle  \\
    & = & \mathcal{L}\left( \langle \sum_{j=0}^n T^jXV^{*j} h , k\rangle \right. \\
    &  & -\left.\langle
     \sum_{j=0}^n T^{j+1}XV^{*j+1} h , k\rangle\right)  \\
     & = & \langle Xh,k\rangle - \mathcal{L}\left( \langle
     T^{n+1}XV^{*n+1}h , k\rangle\right).
\eeq

On the other hand, 
\begin{eqnarray*}
    \langle T^{n+1}XV^{*n+1}h , k\rangle & = & \langle T^{n+1}XV^{*n+2}Vh , k\rangle \\
     & = & \langle \sum_{j=0}^{n+1} T^{j}XV^{*j+1} Vh , k\rangle - \langle 
     \sum_{j=0}^{n} T^{j}XV^{*j+1} Vh , k\rangle .
\end{eqnarray*}
Therefore $\mathcal{L}\left( \langle T^{n+1}XV^{*n+1}h , k\rangle\right) = 0$ 
and thus $X = TZ -ZV$. We also have 
$$\langle ZVV^*h , k\rangle = - \mathcal{L}\left( \langle
\sum_{j=0}^n T^jXV^{*j+1} VV^*h , k\rangle\right)  = \langle Zh , k\rangle ,$$
which completes the proof.
\end{proof}
We obtain the following known result (cf. \cite{foias/williams,wu,conw,C+3}).
\begin{cor}\label{cor:un}
Let $U\in \Bh$ be a unitary, let $T\in \Bk$ be a power bounded operator and let 
$X\in \Bhk$. Then the operator 
$$R(X) =  R(X;T,U) = \left[\begin{array}{cc}
            T & X  \\
            0 & U
        \end{array}\right] \in \Bkoh$$
is similar to $R(0) = T\oplus U$ if and only if $R(X)$ is power-bounded.
\end{cor}
\begin{proof}
If $R(X)$ is similar to $R(0) = T\oplus U$, then $R(X)$ is power-bounded since $R(0)$ is.

Suppose that $R(X)$ is power-bounded. We have 
$$R^n = \left[\begin{array}{cc}
            T^{n} &
\sum_{j=0}^{n-1}T^{j}XU^{n-j-1}  \\
            0 & U^n
        \end{array}\right] 
$$
and thus 
$$ \|\sum_{j=0}^{n}T^{j}XU^{n-j}\| \leq M $$
for a suitable positive constant $M$. This yields 
\beq
\|\sum_{j=0}^n T^jXU^{*j+1}\| & = & \|\sum_{j=0}^n T^jXU^{n-j}U^{*n-j}U^{*j+1}\| \\
            & \leq & \|\sum_{j=0}^{n}T^{j}XU^{n-j}\| \; \|U^{*n+1}\| \\
        & \leq & M .  
\eeq
Using Theorem \ref{thm:bl}, there is a bounded operator $Z$ such that $X = TZ- ZU$. Therefore
$R(X)$ is similar to $R(0) = T\oplus U$.
\end{proof} 
It follows from \cite{cassier} that if $U$ is a coisometry and $T$ is a contraction, 
then $R(X)$ is similar to a
contraction if and only if $R(X)$ is power bounded. However, $R(X)$ is not necessarily similar to 
$R(0)$ as \cite[Example 2.15]{wu} shows.

The following result is a counterpart of Theorem \ref{thm:bl} (see also
\cite[Th. 4.1]{bp}).
Its proof will be omitted.
\begin{thm}\label{thm:blbis}
Let $T\in \Bk$ be a coisometry (i.e. $TT^* = I$), 
let $V\in \Bh$ be a power bounded operator and let 
$X\in \Bhk$.
Then the commutator equation
$$X = TZ -ZV$$
has a bounded solution $Z\in \Bhk$ with $(I - T^*T)Z = 0$ if and only if 
\begin{equation}\label{eq:gr2}
\sup_{n}\|\sum_{j=0}^n T^{*j+1}XV^{j}\| < \infty .
\end{equation}
\end{thm}
The following result is the first characterization of operators $R(X)$ similar to contractions.
\begin{thm}[growth condition]\label{thm:gr}
\conda Let $V\in \Bh$ be an isometry, let $T\in \Bk$ be a right invertible contraction 
and let 
$X\in \Bhk$.
Then the operator 
$$R(X) =  \left[\begin{array}{cc}
            T & X  \\
            0 & V
        \end{array}\right] \in \Bkoh$$
is similar to a contraction if and only if $X = A + F$, with $A, F \in \Bhk$ 
satisfying
$$\sup_{n}\|\sum_{j=0}^n T^jAV^{*j+1}\| < \infty$$
and 
$$FV = 0.$$
\condb Let $T\in \Bk$ be a coisometry, let $V\in \Bh$ be a left invertible contraction 
and let 
$X\in \Bhk$.
Then the operator 
$$R(X) =  \left[\begin{array}{cc}
            T & X  \\
            0 & V
        \end{array}\right] \in \Bkoh$$
is similar to a contraction if and only if $X = A + F$, with $A, F \in \Bhk$ 
satisfying
$$
\sup_{n}\|\sum_{j=0}^n T^{*j+1}AV^{j}\| < \infty 
$$
and 
$$TF = 0.$$
\end{thm}
\begin{proof} We give the proof only for the second part. 
By Foias-Williams' \cite{C+3} theorem, 
$R(X)$ is similar to a contraction if and only if the commutator
equation $X = TZ -ZV$ has a bounded solution $Z\in \Bhk$. 

Suppose the commutator equation
has a bounded solution $Z\in \Bhk$. Set 
$$F = -(I-T^*T)ZV \quad \mbox{ and } \quad A = X - F = TZ - T^*TZV.$$
Then $X= A+F$, $TF = 0$ and 
$A = TD - DV$ with $D = T^*TZ$ satisfying $(I-T^*T)D = 0$. Apply Theorem \ref{thm:blbis}. 

Suppose now that $X$ has a decomposition $X = A+F$ as required. By Theorem \ref{thm:blbis} 
there exists a bounded
operator $D$ such that $A = TD - DV$ and $(I-TT^*)D = 0$. Let $L$ be a left inverse for $V$. Set
$Z = D-FL$. Then
$$
TZ - ZV = TD - DV - TFL + FLV = A + F = X
$$
and the proof is complete.
\end{proof}
We mention also the following result.
\begin{prop}
Let $V\in \Bh$ be a unilateral shift (i.e. a pure isometry). Let $X \in \Bh$ be an operator
such that
$$
\sup_{n}\|\sum_{j=0}^n \left( V^{j+1}XV^{j} - V^{*j}XV^{*j+1}\right) \| < \infty 
$$
Then the commutator equation $X = V^*Z-ZV$ has a bounded solution $Z \in \Bh$. 
\end{prop}
\begin{proof}
Consider $Z \in \Bh$ given by
$$\langle Zh , k\rangle = \mathcal{L}\left( \langle
\sum_{j=0}^n \frac{V^{j+1}XV^{j} - V^{*j}XV^{*j+1}}{2} h , k\rangle\right) .$$
Then 
$$\langle (V^*Z - ZV)h,k\rangle  = \langle Xh,k\rangle - \frac{1}{2}\mathcal{L}\left( 
\langle \left( V^{n+1}XV^{n+1} + V^{*n+1}XV^{*n+1}\right) h , k\rangle \right) .
$$
We have 
\beq
|\langle V^{n+1}XV^{n+1}h , k\rangle| & = & |\langle XV^{n+1}h , V^{*n+1}k\rangle| \\
 & \leq &\|X\|\,\|h\|\,\|V^{*n+1}k\| 
\eeq
and the last term tends to zero 
since $V^{*n}$ tends strongly to $0$ as $n$ tends to $\infty$. We obtain 
that $Z$ is a bounded solution of the commutator equation.
\end{proof}
\section{Zero-product perturbations}
\begin{dfn}
The operator $T\in \Bh$ is said to be a \emph{zero-product perturbation} 
of $C\in \Bh$ by
$E\in \Bh$ if $T = C+E$ and $EC = 0$.
\end{dfn}
The following result shows that the class of 
operators similar to contractions is
stable under zero-product perturbations by operators of spectral radius smaller
than one, in particular by nilpotent operators.
\begin{thm}\label{thm:pert}
A zero-product perturbation
of an operator similar to a contraction by an operator of spectral radius smaller
than $1$ is similar to a contraction.
\end{thm}
\begin{proof}
Let $C\in \Bh$ be an operator similar to a contraction and thus 
completely polynomially
bounded \cite{paulsen:jfa}. This means that there exists a universal constant $K > 0$ such that 
\begin{equation}
 \|P(C)\| \leq K \|P\|_{\infty} = \sup_{|z|=1}\|P(z)\|_{M_{p}(\C)} \label{k}
\end{equation}
for each polynomial $P$ with matrix coefficients in $M_{p}(\C)$, for any $p$.  
Recall that $P(C)$ is identified with an
operator acting on the direct sum of $p$ copies of $\HH$ 
in a natural way.

Let $E \in \Bh$ be an
operator of spectral radius $r(E)$ smaller than $1$ such that
$EC = 0$. Set $T = C+E$.
We have 
$$T^2 = (C+E)^2 = C^2 + CE + E^2$$
and, by recurrence,
\begin{equation}\label{pow}
T^n = C^n + C^{n-1}E + \cdots + CE^{n-1} + E^n ,
\end{equation}
for every $n$. 

Let $p \geq 1$ be fixed. Let 
$$P(z) = \sum_{j=0}^d A_jz^j$$
be a polynomial of degree $d$ 
with $p\times p$ matrix coefficients $A_j \in M_p(\C)$. 

Denote $P_{(0)} =P$, 
$$P_{(1)}(z) = \sum_{j=1}^d A_jz^{j-1} = \frac{P(z)-P(0)}{z}$$
and, recursively,
$$P_{(n)}(z) = \frac{P_{n-1}(z)-P_{n-1}(0)}{z} .$$ 
The representation (\ref{pow}) implies
\begin{equation}
\label{sum} 
P(T) = P(C) + P_{(1)}(C)E + P_{(2)}(C)E^2 + \cdots + P_{(d-1)}(C)E^{d-1}.
\end{equation}

We can estimate the norm 
$$\|P_{(n)}\|_{\infty} = \sup_{|z|=1}\|P_{(n)}(z)\|_{M_p(\C)}$$ 
as follows :
\beq
\|P_{(n)}\|_{\infty} & = & \| z^nP_{(n)}\|_{\infty} \\
        & = & \|P - D_{n-1}\ast P\|_{\infty} .
\eeq
Here 
$$ D_n(t) = \sum_{|j|\leq n} e^{ijt} = 
\frac{\sin ((n+\frac{1}{2})t)}{\sin (t/2)} $$
is the Dirichlet kernel and the convolution $D_{n-1}\ast P$
with the polynomial $P$ with matrix coefficients has an obvious meaning. 

The $L^1$-norm of the Dirichlet kernel grows like $\log n$ \cite[II.(12.1)]{zygmund}.
Therefore there exists a positive constant $A$ such that 
\begin{equation}
\label{a}
\|P_{(n)}\|_{\infty} \leq A \log (n+2) \|P\|_{\infty} , 
\end{equation}
for every $n\ge 0$ and every $p$. 

Combining now equations (\ref{k}), (\ref{sum}) and (\ref{a}), we obtain 
$$
\|P(T)\|_{\infty} \leq AK\left(  \sum_{n=0}^{\infty} \log (n+2) \| E^n\|\right) 
\|P\|_{\infty} .
$$
Since $r(E) < r < 1$, there exists a constant $C$ such that 
$$ \|E^n\| \leq Cr^n$$
for all $n$ 
and thus 
$$\sum_{n=0}^{\infty} \log (n+2) \| E^n\| \leq C\sum_{n=0}^{\infty} r^n\log (n+2)$$
is convergent.
Therefore $T$ is completely polynomially bounded and thus similar to a contraction by 
Paulsen's \cite{paulsen:jfa} criterion. 
\end{proof}
\begin{rmks}
\firsta The series 
$$\sum_{n=0}^{\infty} \log (n+2) \| E^n\|$$
is convergent if and only if $r(E) < 1$. Indeed, 
the convergence of the series with the logarithm implies the convergence of 
$\sum_{n=0}^{\infty}\| E^n\|$ and thus $r(E) < 1$.

\condb For $C = 0$ we obtain the classical theorem of Rota \cite{rota} stating 
that an operator $T$ with $r(T) < 1$  
is similar to a contraction. 

\condc By applying Theorem \ref{thm:pert} to $T^* = C^* + E^*$ we obtain a similar
statement for perturbations satisfying the reversed zero-product condition $CE = 0$.
\end{rmks}
In the sequel we will use only the following corollary for zero-product
perturbations by nilpotents of order two.
\begin{cor}\label{cor:nilp}
Let $C\in \Bh$ be an operator similar to a contraction and let $E \in \Bh$ be a
nilpotent operator such that
$EC = 0$. 
Then $T = C+E$ is similar to a contraction.
\end{cor}
\begin{rmk}
The zero-product 
condition $EC = 0$ (or $CE = 0$) is necessary in Theorem \ref{thm:pert} and in 
Corollary \ref{cor:nilp}. Indeed,
if $\HH$ is the Euclidean space $\C^2$, 
$$ C = I = \left[\begin{array}{cc}
            1 & 0  \\
            0 & 1
        \end{array}\right] \quad \mbox{ and } \quad E = \left[\begin{array}{cc}
            0 & 1  \\
            0 & 0
        \end{array}\right] , $$
then $E^2 = 0$ and $EC = CE = E$. However, 
$$T = C + E = \left[\begin{array}{cc}
            1 & 1  \\
            0 & 1
        \end{array}\right] $$ 
is not similar to a contraction. Indeed, 
$$T^n = \left[\begin{array}{cc}
            1 & n  \\
            0 & 1
        \end{array}\right] $$
is not power bounded. 

Moreover, there exists a contraction $C$ and an operator 
$E$ such that $E^2 = 0$, $EC = CE$, $T = C+E$ is a polynomially 
bounded operator but $T$ is  
not similar to a contraction.
 
For a counterexample, let
$$T = \left[\begin{array}{cc}
            S^{*} & \Gamma  \\
            0 & S^{}
        \end{array}\right] \in 
        \mathcal{B}(\ell^2(\mathcal{H})\oplus \ell^2(\mathcal{H}))$$
be the polynomially bounded operator not similar to a contraction constructed 
by G. Pisier
\cite{pisier:jams}. Here $S$ is the unilateral shift of (infinite) multiplicity $\dim \HH$ 
and $\Gamma$ is a suitable \cite{pisier:jams} 
operator-valued Hankel operator, thus satisfying 
$S^{*} \Gamma = \Gamma S$. 
Then $T$
is the sum of the contraction
$$C = \left[\begin{array}{cc}
            S^{*} & 0  \\
            0 & S
        \end{array}\right]$$
with the operator 
$$E = \left[\begin{array}{cc}
            0 & \Gamma  \\
            0 & 0
        \end{array}\right] $$
satisfying $E^2 = 0$ and $EC = CE$. 
\end{rmk}
\section{Quadratically near operators modulo subspaces}
The following definition was introduced in \cite[Definition 2.5]{b:qn}.
\begin{dfn}\label{dfn:qn}
 Let $\beta : \N \to \R_{+}^{\ast}$. Two operators $T \in \Bh$ and $C \in \Bh$ 
are 
said to be $\beta$-\emph{quadratically near} if 
$$ 
s := \left[\sup_{N\geq 0}\left\|\sum_{n=0}^{N} \frac{1}{\beta(n)^2}(T^n - 
        C^n)(T^n - C^n)^{\ast}\right\|\right]^{1/2}< +\infty .
$$
The two operators are simply called \emph{quadratically near} if 
this condition holds with $\beta(n) = 1$ for each $n$.
We denote $s$ in the above definition by 
$near(T,C,\beta)$. If 
$\beta(n) = 1$ for each $n$, we call $s = near(T,C)$ the \emph{nearness} (or 
the $2$-nearness)
between $T$ and $C$.
\end{dfn}
The following 
result gives equivalent definitions.
\begin{lem}\label{lema}
Let $\beta : \N \to \R_{+}^{\ast}$. 
Two operators $T$ and $C$ in $\Bh$ are 
$\beta$-quadratically near with $near(T,C,\beta) \leq s$ 
\begin{itemize}
\item if and only if 
\begin{equation}\label{eq:21}
\sum_{n=0}^{+\infty} \frac{1}{\beta(n)^2}\left\| (T^n - C^n)^{\ast}y\right\|^2 
\leq s^2 \left\| y\right\|^2  ,
\end{equation}
for all $y \in \HH$ ;
\item if and only if, for every $N \in \N$ and all 
$x_0, \ldots , x_N \in \HH$, we have
\begin{equation}\label{eq:22}
\left\| \sum_{n=0}^{N}\frac{1}{\beta(n)}(T^n - C^n)x_n\right\| \leq s 
 \sqrt{\sum_{n=0}^{N}\left\| x_n\right\|^2}  .
\end{equation}
\end{itemize}
\end{lem}
\begin{proof} The first equivalence was remarked in \cite[Lemma 2.6]{b:qn} and 
follows from the fact that numerical radius of $A$ equals $\|A\|$ for normal 
operators $A$. The second equivalence follows from the fact that Definition \ref{dfn:qn} 
and the condition
(\ref{eq:22}) are both saying that 
$$ \sup_N \|R_N\| \leq s ,$$
where $R_N$ is the row operator
$$
R_N = \left[ 0 \quad \frac{T-C}{\beta(1)} \quad \cdots \quad \frac{T^N-C^N}{\beta(n)} \right] 
$$
acting on column vectors of $\ell^2_N(\HH)$.
\end{proof}
\begin{rmk}
If 
\begin{equation*}
\limsup_{n\to\infty} \| T^n - C^n\|^{1/n} < 1 ,
\end{equation*}
or if 
\begin{equation*}
\sum_{n=0}^{+\infty} \| T^n - C^n\|^2 < +\infty ,
\end{equation*}
then $T$ and $C$ are 
quadratically near. Operators satisfying 
$$\limsup_{n\to\infty} \| T^n - C^n\|^{1/n} = 0$$ 
are called in \cite{apostol2} 
\emph{asymptotically equivalent} operators.
\end{rmk}
It follows from the more general result proved in \cite{b:qn} that if 
$C$ is similar to a contraction and
if $near(T,C) < +\infty$, then $T$ is similar to a contraction. 
In particular, Rota's theorem is
obtained for $C = 0$ since every operator of spectral radius smaller than $1$ is
quadratically near the null operator. We refer to \cite{b:qn} for consequences of the conditon
$near(T,C,\beta) < +\infty$.

We introduce the following definition.
\begin{dfn}
 Let $\HH_0$ be a subspace of $\HH$. Two operators $T$ and $C$ in $\Bh$ 
are said to be $\beta$-\emph{quadratically near modulo} $\HH_0$ if 
for every $N \in \N$ and for all $x_0, \ldots , x_N \in \HH_0$ we have
\begin{equation}\label{eq:mod}
\left\|\sum_{n=0}^{N} \frac{1}{\beta(n)}(T^n - C^n)x_n\right\| \leq s 
\sqrt{\sum_{n=0}^{N}\| x_n\|^2} .
\end{equation}
For $\beta(n) \equiv 1$, we say that $T$ and $C$ are quadratically near modulo $\HH_0$.
\end{dfn}

Quadratically near operators corresponds to $\HH_0 \equiv \HH$. We will indentify $\HH$ with the
subspace $\HH\oplus \{0\}\oplus \{0\}\cdots $ in $\ell^2(\HH)$.

\begin{egz}\label{eg2}
Let $(\omega_k)_{k\ge 0}$ be a sequence of strictly positive weights. 
Denote by $\som$ the unilateral weighted shift operator 
$$\som : \ell^2(\HH) \ni (x_0, x_1 , \ldots) \to (0, \om_0x_0, \om_1x_1 , \ldots) \in
\ell^2(\HH) $$
on $\ell^2(\HH)$. 
Define 
$$\beta(n) = \om_0\cdots \om_{n-1} \mbox{ for } n\ge 1 \mbox{ and }
\beta(0) = 1 .$$ 
Then $\som$ is $\beta$-quadratically near $0$ (the null operator) modulo $\HH$. 
If we have 
$$
\sup_{n\ge 0}\beta(n) < +\infty ,
$$
then $\som$ is quadratically near $0$ 
modulo $\HH$.
\end{egz}
\begin{proof}
Let $x_n  = (x_n, 0, 0, \ldots) \in \HH$, $0\le n \leq N$. Since
$$\sum_{n=0}^N\som^n(x_n, 0, 0, \ldots) = (x_0, \beta(1)x_1, \ldots , \beta(n)x_n , \ldots ),$$
we have 
$$
\|\sum_{n=0}^N\frac{1}{\beta(n)}\som^nx_n\|^2 = 
\sum_{n=0}^N\|x_n\|^2  
$$
and 
$$\|\sum_{n=0}^N\som^nx_n\|^2 = 
\sum_{n=0}^N\beta(n)^2\|x_n\|^2 \le \left[\sup_{k\ge 0}\beta(k)\right]^2 
\sum_{n=0}^N \|x_n\|^2 .$$
\end{proof}

Another example in the same vein is the following.
\begin{egz}\label{eg}
Let $V\in \Bh$ be an isometry. Then $V$ is quadratically near 
$0$ (the null operator) modulo $\ker V^*$.
\end{egz}
\begin{proof} 
According to the Wold decomposition, $V$ is the direct sum of a 
unitary operator $U$ and a pure isometry $S$.
Since $\ker V^* = \{0\}\oplus \ker S^*$, and the unilateral shift 
$S$ is quadratically near $0$ modulo $\ker S^*$ by the preceding example, $U\oplus S$ is 
quadratically near $0$
modulo $\ker V^*$.
\end{proof}

\begin{thm}\label{thm:newom}
Let $(\omega_k)_{k\ge 0}$ be a sequence of strictly positive weights and let $\beta(0) = 1$
and
$$\beta(n) = \om_0\cdots \om_{n-1}   \quad (n\ge 1) .$$
Suppose that
\begin{equation}\label{eq:pwbdd}
\sup_{n,k\ge 0}\frac{\beta(n+k)}{\beta(n)} < \infty .
\end{equation}
Consider $\som \in \B(\ell^2(\HH))$ the unilateral weighted shift with
weights $\om_k$.
Let $T\in \Bk$ be an operator similar to a contraction 
and let 
$X\in \B(\ell^2(\HH),\mathcal{K})$. Suppose the operator 
$$R(X) =  \left[\begin{array}{cc}
            T & X  \\
            0 & \som
        \end{array}\right] \in \B(\KK\oplus \ell^2(\HH))$$
is $\beta$-quadratically near $R(0) = T\oplus \som$ modulo 
$\HH$. 
Then $R(X)$ is similar to a contraction.
\end{thm}
\begin{proof}
The condition (\ref{eq:pwbdd}) means
that $\som$ is power bounded, and thus \cite{shields} $\som$ is similar to a
contraction. Without loss of any generality we can assume that $T$ and $\som$ are two
contractions. In order to prove the similarity of
$R(X)$ to a
contraction we will construct an equivalent, \emph{hilbertian} norm on 
$\B(\KK\oplus\ell^2(\HH))$ such that $R(X)$, with respect to this norm, is a contraction. 
The construction of this new norm is similar to constructions in \cite{holbrook:sz} and
\cite{b:qn}.

Denote 
$$
X_n = \sum_{j=0}^{n-1}T^{j}X\som^{n-j-1} \; .
$$
Since 
$$R(X)^n = \left[\begin{array}{cc}
            T^{n} & X_n  \\
            0 & \som^n
        \end{array}\right] ,
$$
we have, using the $\beta$-quadratic nearness condition,
\begin{equation}\label{eq:star}
\|\sum_{n=0}^N\frac{1}{\beta(n)}X_nu_n\| \le C\sqrt{\sum_{n=0}^N\|u_n\|^2}
\end{equation}
for all $u_n\in \HH$.

Every element $h = (h_0,h_1,\ldots )\in \ell^2(\HH)$ can be (uniquely) written as   
\begin{equation}
h = \sum_{n=0}^{\infty}\frac{1}{\beta(n)}\som^nh_n , \label{dech}
\end{equation}
with $h_n = (h_n,0,0\ldots) \in \HH$.

Consider decompositions of elements $(k,h)$ 
of $\KK\times \ell^2(\HH)$ of the following type~:
decompose $h= (h_0,h_1,\ldots )$ as in (\ref{dech}) ; then decompose $k\in \KK$ as follows~:
\begin{equation}
k = \sum_{n=0}^{\infty}T^nk_n + \sum_{n=0}^{\infty}\frac{1}{\beta(n)}X_nh_n . \label{deck}
\end{equation}

Several remarks are in order about this decomposition.
Firstly, the above series
$$\sum_{n=0}^{\infty}\frac{1}{\beta(n)}X_nh_n$$ 
converges for all $h$. 
Indeed, we have
$$
\|\sum_{n=m}^{m+p}\frac{1}{\beta(n)}X_nh_n\|^2 \le C^2  \sum_{n=m}^{m+p}\|h_n\|^2
$$
and the last sum is bounded by 
the tail of a convergent series. Secondly, we suppose that only a
finite number of the $k_n$'s are non-zero, that is we consider only finite sums in the first
part of the decomposition (\ref{deck}).

Such decompositions of
$(k,h) \in \KK\times \ell^2(\HH)$ always exist. 
Indeed, given the unique decomposition of $h$ as in (\ref{dech}), there is at least one
finite decomposition of $k \in \KK$ as in (\ref{deck})~: take for instance
\begin{equation} 
k_0 = k - \sum_{n=0}^{\infty}\frac{1}{\beta(n)}X_nh_n \mbox{ and } k_n = 0 \mbox{ for } n\ge 1 .
\label{decks}
\end{equation}

We define a new norm $|(\cdot,\ast)|$ on $\KK\times \ell^2(\HH)$ by setting 
\begin{equation}
|(k,h)|^2 = \inf \{ \|\sum_{n\ge 0}T^nk_n\|^2 +
\|h\|^2 + \sum_{n=0}^{\infty}\|k_n\|^2\}\; ,
\end{equation}
where the infimum is taken over all decompositions of $(k,h)$ described above. Note that
$\sum_{n=0}^{\infty}\|k_n\|^2$ is also a finite sum.

We prove now that $|(\cdot,\ast)|$ is a hilbertian norm on 
$\KK\times \ell^2(\HH)$, equivalent to
the usual $\ell^2$ norm on $\KK\oplus\ell^2(\HH)$ and
that $R(X)$ is a contraction with respect to this new norm. 

\begin{gen}$|(\cdot,\ast)|$ {\bf is a seminorm.}
Take two elements $(k,h) , (k',h')$ in $\KK\times \ell^2(\HH)$ with their corresponding 
decompositions as above given by two sequences $(k_n,h_n)$ and  $(k'_n,h'_n)$ in 
$\KK\times \HH$. 
By adding eventually zeros, we may assume that both decompositions have the 
same (finite) number of $k$'s. 
Then $(k+k',h+h')$ is decomposed using the sequence 
$(k'_n +k_n,h'_n + h_n)$, $n\ge 0$.

Using the triangle inequality $\|a+b\| \leq \|a\|+\|b\|$ for 
$$a = (\sum_{n\ge 0}T^nk_n, h, k_0, k_1, \ldots )$$
and
$$b = (\sum_{n\ge 0}T^nk'_n, h', k'_0, k'_1, \ldots )$$ 
and taking the infimum over all representations of $(k,h)$ and $(k',h')$, we get the triangle
inequality for the new norm.

The proofs of the inequality $|\ld (k,h)| \leq |\ld| |(k,h)|$ and its 
converse are left to the reader. 
\end{gen}

\begin{gen}$|(\cdot,\ast)|$ {\bf is an equivalent norm.}
Let $(k,h)\in \KK\times \ell^2(\HH)$ and consider 
a decomposition as above. Then, using (\ref{deck}) and (\ref{eq:star}), we obtain
\beq
\|k\|^2+\|h\|^2 & = & \|\sum_{n\ge 0}T^nk_n +
\sum_{n=0}^{\infty}\frac{1}{\beta(n)}X_nh_n\|^2 + \|h\|^2 \\
 & \leq & 2\|\sum_{n\ge 0}T^nk_n\|^2 + 2\|\sum_{n=0}^{\infty}\frac{1}{\beta(n)}X_nh_n\|^2 +
 \|h\|^2 \\ 
 & \le & 2\|\sum_{n\ge 0}T^nk_n\|^2 + \left( 2C^2+1\right)\|h\|^2 .
\eeq
Taking the infimum over all representations of $(k,h)$, we obtain that 
$\|k\|^2+\|h\|^2$, which is the norm of $(k,h)$ in 
$\KK\oplus\ell^2(\HH)$,  is no greater than a constant times the new norm $|(k,h)|^2$.

For the converse inequality, consider the representation of $h$ as in (\ref{dech}) and of 
$k$ as in (\ref{decks}). Then, using (\ref{eq:star}), we have
\beq
|(k,h)|^2 & \leq & 2\|k-\sum_{n=0}^{\infty}\frac{1}{\beta(n)}X_nh_n\|^2 + \|h\|^2  \\
 & \leq & 4\|k\|^2 +  4\|\sum_{n=0}^{\infty}\frac{1}{\beta(n)}X_nh_n\|^2 + \|h\|^2 \\
& \leq & M\left[\|k\|^2+\|h\|^2\right] ,
\eeq
for a suitable constant $M$.
\end{gen}

\begin{gen}{$|(\cdot,\ast)|$ {\bf is Hilbertian}.} Let
$(k,h) , (k',h') \in \KK\times \ell^2(\HH)$ with their corresponding 
decompositions as above given by $(k_n,h_n), (k'_n,h'_n) \in \KK\times \HH$. 
Then $(k \pm k',h \pm h')$ are decomposed by $(k_n\pm k'_n,h_n\pm h'_n)$. 
Let $\Sigma = |((k + k',h + h')|^2 + |((k - k',h - h')|^2$. Using the parallelogram
law in $\KK$ and $\HH$ we obtain
\beq
\Sigma & \leq & \|\sum_{n\ge 0}T^n(k_n+k'_n)\|^2 + \|h+h'\|^2 + 
\sum_{n=0}^{\infty}\|k_n+k'_n\|^2 \\
 &+& \|\sum_{n\ge 0}T^n(k_n-k'_n)\|^2 + \|h-h'\|^2 +
\sum_{n=0}^{\infty}\|k_n-k'_n\|^2 \\
 & = & 2\left[ \|\sum_{n\ge 0}T^nk_n\|^2 + \|h\|^2 + \sum_{n=0}^{\infty}\|k_n\|^2\right] \\
 &+& 
 2\left[ \|\sum_{n\ge 0}T^nk'_n\|^2 + \|h'\|^2 + \sum_{n=0}^{\infty}\|k'_n\|^2\right] \\
\eeq
Taking the infimum over all representations of $(k,h)$ and $(k',h')$, 
we obtain the parallelogram  
(in)equality for $|(\cdot,\ast)|$. This implies \cite{amir} 
that the norm comes from a scalar product.
\end{gen}

\begin{gen}{\bf The operator} $R(X)$ {\bf with respect to} $|(\cdot,\ast)|$. 
Let $(k,h)\in \KK\times \ell^2(\HH)$ and decompose $k$ and $h$ 
as in (\ref{dech}) and (\ref{deck}). Then 
\beq
\som h & = & (0,\om_0h_0,\om_1h_1, \ldots) \\
& = & \sum_{n\ge 0}\frac{1}{\beta(n+1)}\som^{n+1}\frac{\beta(n+1)}{\beta(n)}h_n 
\eeq
and 
\beq
Tk + Xh & = & \sum_{n\ge 0}T^{n+1}k_n + \sum_{n=0}^{\infty}\frac{1}{\beta(n)}TX_nh_n \\
&+& \sum_{n\ge 0}\frac{1}{\beta(n)}X\som^nh_n\\
& = & \sum_{n\ge 0}T^{n+1}k_n + \sum_{n=0}^{\infty}\frac{1}{\beta(n)}X_{n+1}h_n \\
& = & \sum_{n\ge 0}T^{n+1}k_n +
\sum_{n=0}^{\infty}\frac{1}{\beta(n+1)}X_{n+1}\frac{\beta(n+1)}{\beta(n)} h_n .
\eeq
\end{gen}
Therefore, the second and the first component of  
$$R(X)(k,h) = \left[\begin{array}{cc}
            T & X  \\
            0 & \som
        \end{array}\right]\left[\begin{array}{c}
            k  \\
            h
        \end{array}\right] \in \KK\oplus \ell^2(\HH)$$
are decomposed by 
$(0,\frac{\beta(1)}{\beta(0)}h_0, \frac{\beta(2)}{\beta(1)}h_1, \ldots)$
and, respectively, $(0,k_0, k_1, \ldots)$. Then 
\beq |R(X)(k,h)|^2 & \leq & \|\sum_{n\ge 0}T^{n+1}k_n\|^2 +
\|\som h\|^2 + \sum_{n=0}^{\infty}\|k_n\|^2 \\
& \leq & \|\sum_{n\ge 0}T^{n}k_n\|^2 +
\|h\|^2 + \sum_{n=0}^{\infty}\|k_n\|^2 , 
\eeq
since $T$ and $\som$ are supposed to be contractions. We obtain that 
$R(X)$ is a contraction in the
new norm $|(\cdot,\ast)|$. Therefore $R(X)$ is similar to a contraction.
\end{proof}

\section{Perturbations of $R(0)$}
\begin{thm}[perturbation of $R(0)$]\label{thm:main3}
Let $(\omega_k)_{k\ge 0}$ be a sequence of strictly positive weights and set $\beta(0) = 1$
and $\beta(n) = \om_0\cdots \om_{n-1}$ for $n\ge 1$ .
Consider $\som \in \B(\ell^2(\HH))$ the unilateral weighted shift with
weights $\om_k$.
Let $T\in \Bk$ be a coisometry 
and let 
$X\in \B(\ell^2(\HH),\mathcal{K})$. 
Then 
$$R(X) =  \left[\begin{array}{cc}
            T & X  \\
            0 & \som
        \end{array}\right] \in \B(\KK\oplus \ell^2(\HH))$$
is similar to a contraction if and only if $R(X)$ is power bounded and it is the zero-product
perturbation by a nilpotent (of order two) of an operator $R(A)$ which is 
$\beta$-quadratically near $R(0) = T\oplus \som$ modulo 
$\HH$. 
\end{thm}
\begin{proof}
Suppose that $R(X)$ is similar to a contraction. Then $R(X)$ and $\som$ are power bounded. 
Since $T$
is a coisometry, the Foias-Williams' \cite{C+3} theorem  implies that the commutator
equation $X = TZ -Z\som$ has a bounded solution $Z\in \Bhk$. Let $L$ denote the operator
$$
L : \ell^2(\HH) \ni (y_0,y_1,\ldots) \to 
(\frac{1}{\om_0}y_1, \frac{1}{\om_1}y_2, \ldots ) \in \ell^2(\HH) ,
$$
which is a left inverse of $\som$.

Set 
$$F = TZ(I - \som L) \quad \mbox{ and } \quad A = X - F = TZ\som L - Z\som.$$
Then $F\som = 0$ and $A = TD - D\som$ with $D = Z\som L$. We have 
$$R(X) = \left[\begin{array}{cc}
            T & A  \\
            0 & \som
        \end{array}\right] + \left[\begin{array}{cc}
            0 & F  \\
            0 & 0
        \end{array}\right] =: R(A) + E .$$
The second operator is nilpotent of order two and also
$$ER(A) = \left[\begin{array}{cc}
            0 & F  \\
            0 & 0
        \end{array}\right] \left[\begin{array}{cc}
            T & A  \\
            0 & \som
        \end{array}\right] = 0 .$$
Thus $R(X)$ is a zero-product perturbation by a nilpotent of order
$2$ of 
$$R(A) = \left[\begin{array}{cc}
            T & A  \\
            0 & \som
        \end{array}\right] .$$
        
We show now that $R(A)$ is $\beta$-quadratically near to 
$R(0) = T\oplus \som$ modulo 
$\HH$. Recall that 
$A = TD - DV$ with $D = Z\som L$.
Consider a sequence of elements $x_n = (x_n,0,0\ldots)$ in $\HH$. 
Then $Dx_n = D\som Lx_n = 0$. 
Denote 
$$
A_n = \sum_{j=0}^{n-1}T^{j}A\som^{n-j-1} .
$$
The powers of $R(A)$ are given by 
$$R(A)^n = \left[\begin{array}{cc}
            T^{n} & A_n  \\
            0 & \som^n
        \end{array}\right] .
$$
We obtain 
\beq 
\|\sum_{n=1}^N\frac{R(A)^n - R(0)^n}{\beta(n)} x_n\|  & = &  
\|\sum_{n=1}^N \frac{1}{\beta(n)}A_nx_n\| \\
    & = & \|\sum_{n=1}^N \frac{1}{\beta(n)}\sum_{j=0}^{n-1}T^{j}(TD-D\som)\som^{n-j-1}x_n\| \\
    & =& \|\sum_{n=1}^N\frac{1}{\beta(n)}\left( T^nD - D\som^n \right) x_n \| \\
    & = & \| - D\sum_{n=1}^N \frac{1}{\beta(n)}\som^nx_n\|\\
    & \leq & \|D\| \left( \sum_{n=1}^N \|x_n\|^2\right)^{1/2} .
\eeq
In the last line we have used Example \ref{eg2}. 
Thus $R(A)$ is $\beta$-quadratically near $R(0)$ modulo $\HH$.

For the converse implication, suppose that $R(X) = R(A) + E$ can be decomposed as in the theorem. 
We obtain $X = A + F$ with $F\som = 0$. 
Note that the power boundedness of 
$R(X)$ implies that of $\som$. By Theorem \ref{thm:newom}, the operator
$$R(A) = \left[\begin{array}{cc}
            T & A  \\
            0 & \som
        \end{array}\right] \in \Bkoh$$ 
is similar to a contraction. By Corollary \ref{cor:nilp}, the zero-product perturbation $R(A) + E$ is also
similar to a contraction.
\end{proof}
\begin{thm}[perturbation of $R(0)$ again]
\label{thm:main1}
Let $V\in \Bh$ be an isometry, let $T\in \Bk$ be a right invertible contraction 
and let 
$X\in \Bhk$. Then the operator 
$$R(X) =  \left[\begin{array}{cc}
            T & X  \\
            0 & V
        \end{array}\right] \in \Bkoh$$
is similar to a contraction if and only if $R(X)$ is power bounded and it is the zero-product 
perturbation by
a nilpotent of an operator
$R(A)$ which is quadratically near $R(0) = T\oplus V$ modulo 
$\ker V^*$.
\end{thm}
\begin{proof}
Suppose that $R(X)$ is similar to a contraction. Since $V$ is an isometry, the commutator
equation $X = TZ -ZV$ has \cite{C+3} a bounded solution $Z\in \Bhk$. Set 
$$F = TZ(I - VV^*) \quad \mbox{ and } \quad A = X - F = TZVV^* - ZV.$$
Then $FV = 0$ and $A = TD - DV$ with $D = ZVV^*$. The proof that the decomposition
$$R(X) = \left[\begin{array}{cc}
            T & A  \\
            0 & V
        \end{array}\right] + \left[\begin{array}{cc}
            0 & F  \\
            0 & 0
        \end{array}\right] $$
fulfils all the requirements is similar to that given above.

For the converse implication, suppose that $R(X) = R(A) + E$ can be decomposed as the theorem. 
In particular $E^2 = 0$ and $X = A + F$ with $FV = 0$. By Corollary \ref{cor:nilp} 
it is sufficient to
show that 
$$R(A) = \left[\begin{array}{cc}
            T & A  \\
            0 & V
        \end{array}\right] \in \Bkoh$$ 
is similar to a contraction. 

Now $R(A) = R(X) - E$ and $ER(X) = 0$ since $FV = 0$. 
This implies that $R(A)^n = R(X)^n - R(X)^{n-1}E$ and thus $R(A)$ is power bounded since $R(X)$ is.

Writing the Wold decomposition for $V$, we get the orthogonal sum $\HH =
\HH_s\oplus \HH_u$ in which $\HH_s$ and $\HH_u$ reduce $V$. Denote by $U$ 
the part of $V$ on $\HH_u$
($U$ is unitary) and by $S$ the part of $V$ on $\HH_s$ which 
is a pure isometry (a unilateral shift). 

With respect to the decomposition $\mathcal{K}\oplus\HH =
\mathcal{K}\oplus\HH_s\oplus\HH_u$, the matrix of $R(A)$ is given by
$$R(A) =  \left[\begin{array}{ccc}
            T & A_s & A_u \\
            0 & S & 0 \\
            0 & 0 & U
        \end{array}\right]  .$$ 
Recall that $U$ is unitary and $R(A)$ is power bounded. By Corollary \ref{cor:un}, 
$R(A)$ is similar to $C\oplus U$, where 
$$C =  \left[\begin{array}{cc}
            T & X_s  \\
            0 & S
        \end{array}\right] \in \mathcal{B}(\mathcal{K}\oplus\HH_s) .$$
Now $R(A)$ quadratically near $R(0)$ modulo $\ker V^*$ implies that $C$ is quadratically near
$T\oplus S$ modulo $\ker S^*$ and thus $C$ is similar to a contraction by Theorem 
\ref{thm:newom}.
The proof is complete.
\end{proof}
\begin{rmk}     
Another proof of Theorem \ref{thm:main1} can be obtained by using a result due to
C. Apostol \cite{apostol:pams}. With our terminology, Apostol proved that if 
$T$ and $X$ are arbitrary and $V$ is a unilateral shift, then the commutator 
equation $X = TZ -ZV$ has
a bounded solution $Z$ satisfying $Z(I-VV^*) = 0$ if and only if $R(X)$ is quadratically
near $R(0)$ modulo $\ker V^*$. Apostol's result, 
together with Theorem \ref{thm:bl}, shows that 
if $V$ is a unilateral shift and $T$ is power bounded, then  $R(X)$ is quadratically
near $R(0)$ modulo $\ker V^*$ if and only if the growth condition (\ref{eq:gr}) holds.
It seems that Apostol's proof does not generalize to
weighted shifts and this explains our proofs of 
Theorems \ref{thm:newom} and \ref{thm:main3}.
\end{rmk}

The conclusion of Theorem \ref{thm:main1} can be strengthened if the 
isometry $V$ is
supposed to be a unilateral shift. Indeed, combining the result of 
C. Foias and J.P. Williams \cite{foias/williams}, the
result of Apostol \cite{apostol:pams}, and our previous results we obtain the
following characterization.
\begin{thm}[nearness plus admissible perturbations]
\label{thm:pure}
Let $V\in \Bh$ be a pure isometry (a unilateral shift), 
let $T\in \Bk$ be a right invertible contraction 
and let 
$X\in \Bhk$. Then the operator 
$$R(X) =  \left[\begin{array}{cc}
            T & X  \\
            0 & V
        \end{array}\right] \in \Bkoh$$
is similar to a contraction if and only if 
$R(X)$ is the zero-product 
perturbation by
a nilpotent of an operator
$R(A)$ which is quadratically near $R(0) = T\oplus V$ modulo 
$\ker V^*$.
\end{thm}

The difference between this characterization and Theorem \ref{thm:main1} is that the
condition of power-boundedness is now missing. 

\section{Concluding Remark}
Sometimes the nearness conditions suffice in the
characterization given by Theorem \ref{thm:pure}. We will show that this happens for the class of
operators studied by Pisier \cite{pisier:jams} and Davidson-Paulsen \cite{DaPa}. 

Let
$\Lambda$ be a function from $\HH$ into
$\Bh$ satisfying
the CAR - \textit{canonical anticommutation relations}~:
 for all $u,v\in H$,
\begin{gather*}\label{CAR}
 \Lambda(u)\Lambda(v)+\Lambda(v)\Lambda(u)=0  \\
\intertext{and}
 \Lambda(u)\Lambda(v)^*+\Lambda(v)^*\Lambda(u)= (u,v) I .
\end{gather*}
The range  of $\Lambda$ is isometric to Hilbert space.
Let $\{e_{n}\}_{n\geq 0}$ be an orthonormal basis for $\HH$, and let
$C_n=\Lambda(e_n)$ for $n\geq 0$. For an
arbitrary sequence $\alpha=(\alpha_0,\alpha_1,\dots)$ in $\ell^2$, let
$$
 \Gamma_\alpha = \Bigl[ \alpha_{i+j}C_{i+j} \Bigr]
$$
be a CAR-valued Hankel operator. Let
$$
R(\Gamma_{\alpha}) = R(S^{*},S ; \Gamma_\alpha) =
   \begin{bmatrix}S^{*}&\Gamma_\alpha\\0&S\end{bmatrix} \in
   \mathcal{B}(\ell^2(\mathcal{H})\oplus \ell^2(\mathcal{H}))
$$
be the
corresponding CAR-valued Foguel-Hankel operator \cite{pisier:jams,DaPa}. 
Here $S \in \mathcal{B}(\ell^2(\mathcal{H}))$ denotes the unilateral forward shift of multiplicity 
$\dim \HH$.

For a fixed sequence $\alpha=(\alpha_0,\alpha_1,\dots) \in \ell^{2}$,
let
$$A(\alpha) = \sup_{k\ge0} (k+1)^2 \sum_{i\ge k}|\alpha_i|^2$$
and
$$B(\alpha) = \sum_{k\ge0} (k+1)^2 |\alpha_k|^2 .$$
The operator $R(\Gamma_{\alpha})$ is \cite{pisier:jams,DaPa} polynomially bounded 
if and only if
$A(\alpha)$
is finite. On the other hand, $R(\Gamma_{\alpha})$ is \cite{pisier:jams,DaPa,ri} similar
to a contraction if and only if 
$B(\alpha)$ is finite (see also \cite{b:qn,bp}).

\begin{thm}
\label{thm:suffice}
Let $\alpha=(\alpha_0,\alpha_1,\dots) \in \ell^{2}$. The operator $R(\Gamma_{\alpha})$ is similar
to a
contraction if and only if it is quadratically near $R(0) = S^{*}\oplus S$ modulo 
$\mathcal{H}$.
\end{thm}
\begin{proof}
Using Theorem \ref{thm:newom}, the nearness condition implies similarity to a contraction. 
For the converse implication, note that
$S^{*}\Gamma_{\alpha} = \Gamma_{\alpha}S$. This implies that 
$$R(\Gamma_{\alpha})^n =  \left[\begin{array}{cc}
            S^{*n} & \Gamma_n  \\
            0 & S^{n}
        \end{array}\right] ,$$
where
$$ \Gamma_n = \sum_{j=0}^{n-1}S^{*j}\Gamma_{\alpha}S^{n-j-1} =
n\Gamma_{\alpha}S^{n-1} .$$

We have 
\beq 
\sum_{n=1}^N \left[ R(\Gamma_{\alpha})^n - R(0)^n\right](h_n, 0,0, \cdots) & = & 
\sum_{n=1}^Nn\Gamma_{\alpha}
S^{n-1}(h_n, 0,0, \cdots)\\
 & = & \Bigl[ (j+1)\alpha_{i+j}C_{i+j} \Bigr]_{i,j\geq 0}\left[\begin{array}{c}
            h_1  \\
            \vdots\\
            h_N
        \end{array}\right]
\eeq

If the operator $R(\Gamma_{\alpha})$ is similar
to a contraction, then \cite{DaPa} $B(\alpha)$ is finite, and thus \cite[Proposition 1]{ri} 
the norm of the matrix 
$$ \Bigl[ (j+1)\alpha_{i+j}C_{i+j} \Bigr]_{i,j\geq 0}$$
is bounded by $B(\alpha)^{1/2}$. We obtain 
$$\| \sum_{n=1}^N \left[ R(\Gamma_{\alpha})^n - R(0)^n\right]x_n\| 
\leq B(\alpha)^{1/2}\sqrt{\sum_{n=1}^{N}\| x_n\|^2}$$
for all $x_n\in \HH$. Therefore, $R(\Gamma_{\alpha})$ is quadratically 
near $R(0) = S^{*}\oplus S$ modulo $\mathcal{H}$.
\end{proof}

\end{document}